\documentclass[11pt]{amsart}

\usepackage{graphicx, amsmath}

\newtheorem{thm}{Theorem}[section]
\newtheorem{prop}[thm]{Proposition}
\newtheorem{lma}[thm]{Lemma}

\theoremstyle{definition}
\newtheorem{definition}[thm]{Definition}
\newtheorem{example}[thm]{Example}

\newcommand{\R}{\mathbf{R}}
\renewcommand{\S}{\mathbf{S}}

\begin{document}
\parindent 0pt
\parskip 5pt
\baselineskip 14pt

\title{\textbf{Projecting $(n-1)$-cycles to
zero on hyperplanes in $\R^{n+1}$}}

\author{Bruce Solomon}
\address{Department of Mathematics, Indiana University,
Bloomington, IN 47405}
\email{solomon@indiana.edu}
\urladdr{php.indiana.edu/$\sim$solomon}

\date{{\tiny July 2001}}



\def\a{\alpha}
\def\g{\gamma}
\def\gdot{\dot\gamma}
\def\G{\Gamma}
\def\O{\Omega}
\def\s{\sigma}
\def\<{\langle}
\def\>{\rangle}
\def\Mbar{\overline M}
\def\D{{\mathcal{D}}}
\def\O{\Omega}
\def\cur#1{{\bf [}#1{\bf ]}}
\def\mass{{\bf M}}
\def\norm{{\bf N}}
\def\spt{\mathrm{spt}}
\def\bd{\partial}
\def\psharp#1{({\bf p}_{#1})_\#}
\def\eps{\varepsilon}
\def\rp{\rho_{_P}}
\def\dbl#1{\mathrm{dbl}_{#1}}

\begin{abstract}
\noindent The projection of a compact oriented submanifold
$\,M^{n-1}\subset \R^{n+1}\,$ on a hyperplane $\,P^{n}\,$ can
fail to bound any region in $\,P\,$. We call this ``projecting
to zero.'' Example: The equatorial $\,\S^1\subset
\S^{2}\subset\R^3\,$ projects to zero in any plane containing
the $\,x_{3}$-axis. Using currents to make this
precise, we show: \emph{A lipschitz (homology) $(n-1)$-sphere
embedded in a compact, strictly convex hypersurface cannot
project to zero on $\,n+1\,$ linearly independent hyperplanes
in $\,\R^{n+1}\,$.} We also show, using examples, that all the
hypotheses in this statement are sharp.
\end{abstract}

\maketitle

\section{Introduction.}
\label{sec:intro}

Basic differential topology shows that a smooth compact
submanifold $\,M^{n-1}\,$ embedded in $\,\R^{n}\,$ always
bounds a domain. But when we embed $\,M\,$ in
$\,\R^{n+1}\,$, and then project it orthogonally into a
hyperplane $\,P^{n}\,$,
\[
M^{n-1}\subset
\R^{n+1}\stackrel{\pi}{\longrightarrow} P^{n}
\]
the projection $\,\pi(M)\,$ will typically bound a linear
combination of simple domains in $\,P\,$ with ``winding
number'' coefficients.  In certain non-generic circumstances
all of these winding numbers can vanish, and in such
cases, we will say that $\,M\,$ \emph{projects to zero on
$\,P\,$}.

For instance, if we embed the round circle $\,\S^1\,$ into a
horizontal plane in $\,\R^3\,$, it projects to zero on any
vertical 2-plane; the projection ``cancels itself''  by
traversing a single line segment once in each direction
(figure 1a).

Figure 1b depicts a subtler---perhaps even
surprising---way in which cancellation can occur. There, all
winding numbers vanish, even though the projected loop is
\emph{immersed}---a pinched figure-8 with each lobe
traversed once in both directions.

\begin{figure}
\includegraphics{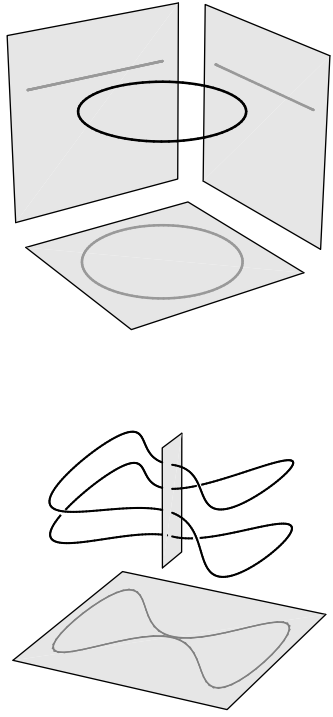} 
\caption{(a) A horizontal circle projects to zero on vertical
(but not horizontal) planes. (b) The projection of this embedded
loop on the horizontal plane is immersed, yet vanishes.}
\end{figure}

While examples like these are not hard to construct,  we
direct our efforts here toward the difficulties that arise
when one tries to make a compact embedded $(n-1)$-manifold
project to zero \emph{simultaneously on a maximal collection
of independent hyperplanes in $\,\R^{n+1}\,$.} We call a
set of hyperplanes \emph{independent} if their normals form a
linearly independent set.

Our main result, Theorem~\ref{thm:main} below, isolates sharp
conditions that obstruct this kind of simultaneous
null-projection. In its statement, a \emph{strict $\,C^2\,$
ovaloid} means a compact, convex $\,C^2\,$ hypersurface
with no vanishing principal curvatures.

{\bf Theorem \ref{thm:main}} \emph{A lipschitz embedding of a
homology $(n-1)$-sphere on a strict $\,C^2\,$ ovaloid in
$\,\R^{n+1}\,$ cannot project to zero on $\,n+1\,$
independent hyperplanes.}

To prove this, we need a more precise definition of
``projection to zero,'' and in Section \ref{sec:currents}, we
employ the language of \emph{currents} for that
purpose. For now, however, we forego rigor, and give
some simple geometric examples that illustrate the problem,
and justify the hypotheses in our theorem above.


\section{Examples}
\label{sec:eg}

The first three examples below highlight topological issues
that arise in the problem we investigate.

\begin{example}
\label{eg:connected}
(Connectedness) Failing connectedness, we
can make a compact embedded hypersurface project to zero on
\emph{any} finite collection of hyperplanes, no matter how
numerous.

Indeed, suppose we have $\,N\,$ hyperplanes
$\,P_1,\,P_2,\,\dots,\,P_N\,$ in $\,\R^{n+1}\,$, with
corresponding unit normals
$\,\nu_1,\,\nu_2,\,\dots,\,\nu_N\,$. Take any embedded,
oriented submanifold $\,M_0\,$, and recursively define, for
$\,i=1,\,2\,\dots,\,N\,$,
\[
M_{i} = M_{i-1}\ \cup\ \left(-M_{i-1} + c_i\,\nu_{i}\right)\
,
\]
Here $\,-M_{i-1}+ c_i\,\nu_i\,$ denotes an
orientation-reversed copy of $\,M_{i-1}\,$, translated
by $\,c_i\,$ units in the $\,\nu_i\,$ direction. By choosing
each $\,c_i\,$ large enough to make the union disjoint, we can
preserve embeddedness throughout the construction.

Projection onto $\,P_{i}\,$ annihilates $\,\nu_{i}\,$, so
$\,M_{i}\,$ projects to zero on $\,P_{i}\,$. Neither
orientation-reversal, nor translation by $\,\nu_{i+1}\,$  in
the next step of the construction will destroy this property,
so after $\,N\,$ iterations, we obtain an embedded submanifold
$\,M_N\,$ that projects to zero on all the $\,P_i$'s, as
claimed.

Figure 2 illustrates this construction with $\,n=0\,$. The
signed circles locate oriented points ($0$-currents) in
the plane. Higher dimensions (or even higher codimensions)
present no additional complications, and the construction
clearly shows that without connectedness, nothing
obstructs simultaneous projection to zero on any finite
number of hyperplanes.

\begin{figure}
\includegraphics{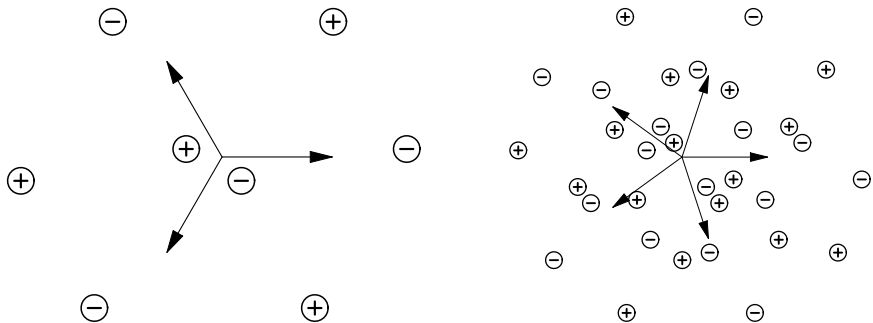} 
\caption{(a) A collection of $\,2^3\,$ signed points that
projects to zero in 3 directions. (b) A collection of $\,2^5\,$
points projecting to zero in 5 directions.}
\end{figure}

\end{example}
\medskip

\begin{example}
\label{eg:embed}
(Embeddedness) Failing embeddedness, one can ``fake''
connectedness to exploit the same phenomenon just discussed.
For instance, when $\,n=2\,$, the latitudinal circles $\,x_3 =
\pm 1/2\,$ in $\,S^2\,$, if oppositely oriented, satisfy all
the hypotheses of Theorem \ref{thm:main} except the
homological one. They violate the conclusion by
projecting to zero on all $\,n+1 = 3\,$ coordinate planes.

The homological defect can be fixed at the expense of
embeddedness, however: Connect these circles with a doubled
longitudinal arc, parametrized once from the lower circle to
the upper, and once from the upper to the lower. The doubled
arc already cancels itself out in $\,\R^3\,$, so the same
holds for its projections. But the resulting connected
curve---a lipschitz immersion of $\,\S^1\,$ into the
sphere---still projects to zero on all three
coordinate planes (figure 3b).
\end{example}

\begin{figure}
\includegraphics{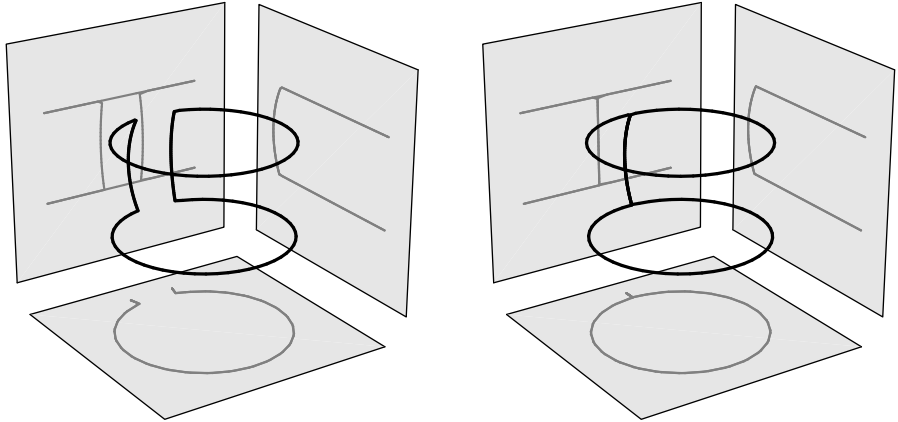} 
\caption{(a) An embedded circle that lies on $\,\S^2\,$ and
projects to zero on 2 coordinate planes. (b) Closing the gap
between the two longitudinal arcs, we make the loop project to
zero on 3 independent planes, but lose embeddedness. (The
horizontal projection vanishes by tracing the circle once each
way.)}
\end{figure}

\begin{example}
\label{eg:torus}
(Topological type) Theorem \ref{thm:main} specifies all the
homology groups of $\,M\,$. To see that we need some such
hypothesis beyond simply $\,H_0(M)= \{0\}\,$ (i.e.,
connectedness), consider the Clifford torus
\[
T^2:=\S^1\times\S^1\subset\R^2\times\R^2\approx\R^4
\]
This torus clearly lies in the origin-centered sphere of
radius $\sqrt{2}$ in $\,\R^4\,$, which is certainly a strict
$\,C^2\,$ ovaloid. Parametrizing $\,T^2\,$ by the map
\[
X(u,v) = \left(\cos u,\ \sin u,\ \cos v,\ \sin v\right)\ ,
\]
one immediately sees that its projection onto any of the
$\,n+1=4\,$ coordinate 3-planes constitutes a cylinder
$\,S^1\times [-1,1]\,$. But in each case, the parametrization
 traverses the interval $\,[-1,1]\,$ once in each direction,
doubling the cylinder with opposite orientations. So all four
coordinate projections vanish as currents.
\end{example}

\begin{example}
\label{eg:dependence}
(Linear independence)  The equator $\,x_{n+1}\equiv 0\,$ in
$\,S^{n}\,$ clearly projects to zero on all of the infinitely
many hyperplanes that contain the $\,x_{n+1}$-axis. Since one
can make $\,n\,$ of these---but no more---linearly
independent, our theorem cannot specify fewer projections
than it does. Figure 1a depicts the case $\,n=2\,$.
\end{example}

\begin{example}
\label{eg:strict}
(Strict convexity) Finally, and perhaps most suprisingly, one
cannot omit the \emph{strict} convexity assumption;
convexity alone does not suffice. Figure 4 presents two
loops---depicted as boundaries of polyhedral surfaces on the
unit  cube---which project to zero on all three standard
coordinate planes.

\begin{figure}
\label{eg:snakes}
\includegraphics{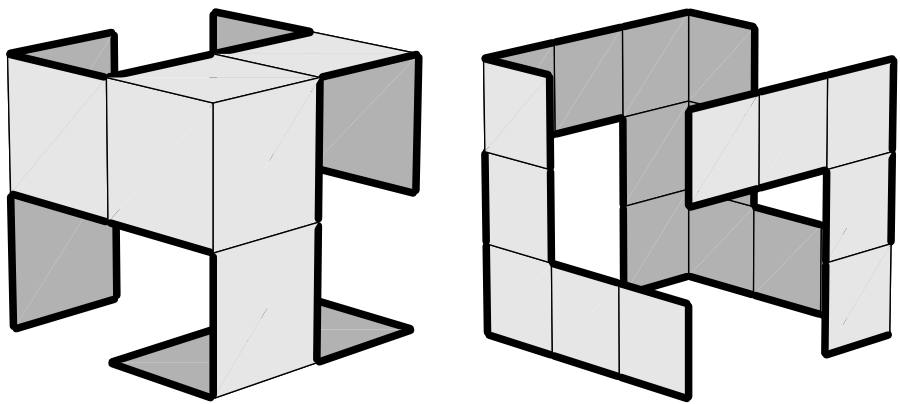} 
\caption{The boundaries of these two polyhedra satisfy all
hypotheses of Theorem~\ref{thm:main} except for \emph{strict}
convexity of the surface---a cube---on which they lie. Both
boundary loops project to zero on all 3 coordinate planes.}
\end{figure}

\end{example}

The loops in Figure 4 have corners, and one cannot simply
round them off without destroying at least one of the null
projections that make them interesting. For this reason, we
originally guessed that no \emph{smooth} loop could project
to zero in three independent directions, even without the
strict convexity assumption in Theorem~\ref{thm:main}.
Mohammad Ghomi has disproven this conjecture, however, by
constructing an intricate counterexample on a ``rounded''
cube. His example appears as an appendix \cite{ghomi:smooth}
to the present paper.
\bigskip

\section{Lipschitz curves and Currents.}
\label{sec:currents}

\begin{definition}
\label{df:lipcurve}
By an \emph{oriented lipschitz $k$-chain} in a riemannian
manifold $\,N\,$, we mean a lipschitz mapping $\,F:M^k\to
N\,$, in which $\,M\,$ is an oriented $k$-dimensional
riemannian manifold. We call $\,F\,$ a (compact)
\emph{$k$-cycle} when $\,M\,$ is closed (compact).
\end{definition}

Our main interest here lies with compact lipschitz
$(n-1)$-cycles in $\,\R^{n+1}\,$ that \emph{project to zero}
on several hyperplanes, as discussed in the previous two
sections. We now formulate that notion more precisely by
viewing lipschitz $k$-chains as $k$-dimensional
\emph{currents}:
\medskip

\begin{definition}
\label{df:cur} (Currents)
Suppose $\,N\,$ is a riemannian manifold. For the purposes of
this paper, a \emph{$k$-current} in $\,N\,$ is simply a
bounded linear function on the vector space $\,\D^k(N)\,$
comprising all compactly supported smooth differential
$\,k$-forms on $\,N\,$.

This definition encompasses some very general objects, but the
currents that interest us here all arise from oriented
lipschitz $k$-chains in a very simple way. Indeed, any
such chain $\,F:M\to N\,$ induces a $k$-current
$\,[F]\,$, via
\[
[F](\phi)
:= \int_{{M}}F^*\phi
 \qquad\hbox{for all }\phi\in\D^k(N)\ ,
\]
Since lipschitz mappings are differentiable almost everywhere
(Rademacher's Theorem (\cite[3.1.6]{federer:GMT}), the
integral makes sense. Further, the change-of-variable formula
\cite[3.2.6]{federer:GMT}, shows that composing $\,F\,$ with
any orientation-preserving diffeomorphism of $\,M\,$ (i.e.,
``reparametrizing'') leaves $\,[F]\,$ unchanged. In this
sense, the current $\,[F]\,$ is actually a more ``geometric''
object than the mapping $\,F\,$.
\end{definition}

Currents also enjoy an elegant notion of boundary:

\begin{definition}
\label{df:bd} (Boundaries)
The \emph{boundary} of a $k$-current $\,T\,$ in a manifold
$\,N\,$ is the $(k-1)$-current $\,\bd T\,$ characterized by
\[
\bd T(\phi) := T(d\phi)\ .
\]
This is an very natural and geometric definition, because
Stokes' theorem combines with Definition \ref{df:cur}
to show that whenever $\,F:M\to N\,$ is an immersed oriented
submanifold with boundary, we have
\[
\bd\cur{F} = \cur{F|_{\bd M}}\ .
\]
\end{definition}
\medskip

\begin{definition}
\label{df:push} (Maps of currents)
Given any  $k$-current $\,T\,$ in $\,N\,$, and a locally
lipschitz mapping $\,G:N\to N'$ between riemannian manifolds,
we get a new current $\,G_\#T\,$ in $\,N'$ via
\[
G_\#T := T\circ G^*
\]
Like the boundary operator, this notion of mapping is
geometrically natural, in the sense that when $\,F:M\to N\,$
is an oriented lipschitz chain, we have $\,G_\#[F] =
[F\circ G]\,$. Moreover, lipschitz mapping commutes with the
boundary operator for currents, just as diffeomorphisms
preserve boundaries of manifolds. More precisely, whenever
$\,G:M\to N\,$ is a proper lipschitz map, and $\,T\,$ is a
$k$-current in $\,M\,$, we have
\begin{equation}
\label{eqn:commute}
\bd\left(G_\#T\right) = G_\#\left(\bd T\right)\ .
\end{equation}
This fact is proven and discussed in
\cite[4.1.14]{federer:GMT}; we shall need it more than once.
\end{definition}
\medskip

Finally, we consider the most basic relationship between
currents and point-sets. The point-set associated with a
mapping is its \emph{image}. The analogous set associated
with a current is its \emph{support}:

\begin{definition}
\label{df:spt} (Supports) Suppose $\,T\,$ is a current on
a manifold $\,N\,$. We say that $\,T\,$ \emph{vanishes} (or
$\,T\equiv 0\,$) on an open subset $\,U\subset N\,$ if
$\,T(\phi)=0\,$ for all 1-forms $\,\phi\,$ with support
$\,\spt(\phi)\subset U\,$. We then define $\,\spt(T)\,$, the
\emph{support} of $\,T\,$, as the complement of the largest
such open set. Equivalently,
$$
\spt(T) := \left\{\,p\in N\ :\  T\not\equiv 0\ \text{on every
neighborhood of $\,p\,$}\,\right\}
$$
Note that the support of a current is always closed.
\end{definition}

We can now give precise
meaning to our concept of ``projection to zero'':

\begin{definition} (Projection to zero)
\label{df:2zero}
When  $\,\pi:\R^{n+1}\to P\,$ is orthogonal projection onto an
affine subspace $\,P\,$, and $\,T\,$ is a current in
$\,\R^{n+1}\,$ such that $\,\pi_\#T\equiv 0\,$, we
say that \emph{$\,T\,$ projects to zero on $\,P\,$.}
\end{definition}

\begin{example}
Consider any smooth oriented $k$-cycle $\,F:M^k\to
\R^{n+1}\,$. As geometric intuition suggests, we get zero
when we project the associated $k$-current $\,[F]\,$ into a
subspace $\,P\,$ of dimension $k$ or less. For, every
$k$-form $\,\phi\,$ on such a subspace is exact, making the
pull-back $\,\pi^*\phi\,$ likewise exact. The projected
current $\,\pi_\#[F] = [\pi\circ F]\,$ then vanishes  because
an exact $k$-form on a compact $k$-manifold always integrates
to zero.
\end{example}

Examples like this show that the support of
a current $\,[\pi\circ F]\,$ can be empty, while the
{image} of the mapping $\,\pi\circ F(M)\,$ is simultaneously
large. In particular, the support of the current induced by a
$k$-chain certainly need not \emph{coincide} with the image of
that $k$-chain. But the two sets  \emph{will} coincide when
$\,M\,$ is compact and $\,F:M\to N\,$ is \emph{injective}.
This fact plays a key role in our work, because it connects
an analytic property of the current $\,[F]\,$ to the topology
of $\,F\,$ itself. We state and prove it as follows:

\begin{lma}
\label{lma:spt}
When $\,M\,$ is compact and $\,F:M\to N\,$ is an
\emph{injective} lipschitz $k$-chain, we have
$\,\mathrm{\spt}([F])=F(M)\,$.
\end{lma}

\begin{proof}
When $\,M\,$ is compact, $\,F(M)\,$ is closed. So any point
$\,p\not\in F(M)\,$ has a neighborhood separating it from
$\,F(M)\,$, and, clearly, excluding $\,p\,$
from $\,\spt([F])\,$. This implies half of our lemma:
\emph{When $\,M\,$ is compact and $\,F:M\to N\,$ is a
lipschitz $k$-chain, we have $\,\spt([F])\subset F(M)\,$}.

Though injectivity was not used above, it is crucial for
the reverse inclusion $\,F(M)\subset\spt([F])\,$. The latter
holds trivially when $\,k=0\,$, so we proceed by induction:
Assuming the inclusion in dimension $\,k-1\,$, we argue that
it also holds in dimension $\,k\,$.

Suppose, toward a contradiction, that it failed for some
$k$-chain $\,F:M^k\to\R^{n+1}\,$. Then for some $\,x\in M\,$,
$\,F(x)\not\in\spt([F])\,$. Since $\,\spt([F])\,$ is closed
and $\,F\,$ is continuous, we then get a non-empty open ball
$\,B\subset M\,$ with
\begin{equation}
\label{eq:spt}
x\in B\quad\text{and}\quad
F(\overline B)\cap \spt([F])=\emptyset\ .
\end{equation}
On the one hand, this implies that
\begin{equation}
\label{eq:vanishing}
[F|_B] = 0\ .
\end{equation}
Indeed, since $\,M\,$ is compact, our injectivity assumption
makes $\,F\,$ an open map. So for some open
$\,U\subset N\,$, we have
\[
F(M)\cap U = F(B)\ .
\]
If $\,\cur{F|_B}\,$ didn't vanish, this would imply that for
some $k$-form $\,\phi\,$ supported in $\,U\,$, we have
$\,
\cur{F}(\phi) = \cur{F|_B}(\phi)\ne 0
\,$,
contradicting Equation (\ref{eq:spt}) above.

On the other hand, once Equation (\ref{eq:vanishing}) holds,
we can deduce $\,\cur{F|_{\bd B}} = \bd\cur{F|_B} = 0\,$,
because the boundary operator commutes with proper lipschitz
mappings (Equation (\ref{eqn:commute}) above). It then
follows that $\,\spt\left(\cur{F|_{\bd B}}\right)=
\emptyset\not\supset F(\bd B)\,$.  Since $\,\bd B\,$ is a
$(k-1)$-dimensional sphere, and $\,F|_{\bd B}:\bd B\to N\,$
is a lipschitz injection, this contradicts our induction
hypothesis.
\end{proof}

\section{Ovaloids.}
\label{sec:ovaloids}

Our main theorem governs lipschitz cycles on
{ovaloids}.

\begin{definition}
\label{df:ovaloid}
A $\,C^2\,$ hypersurface $\,Q\subset {\R^{n+1}}\,$  is an
\textit{ovaloid} if it bounds a compact, convex domain. We
call an ovaloid \textit{strict} if is strictly convex, i.e.,
when its outward unit normal (gauss) mapping
\[
\nu:Q\to{\S^n}
\]
is a diffeomorphism. Equivalently, $\,Q\,$ is strict when
it has everywhere positive principal curvatures.
\end{definition}

When an ovaloid $\,Q\,$ is symmetric with respect to
reflection across a hyperplane $\,P\subset\R^{n+1}\,$, the
symmetry induces a smooth involution $\,\rho:Q\to Q$.
Choosing a unit vector $\, u\,$ normal to $\,P$, we can
express this symmetry by the formula
\begin{equation}
\label{eq:r1}
\rho(x) = x - 2(x\cdot u)\, u \ .
\end{equation}
In this situation, the symmetry hyperplane meets $\,Q\,$ along
the zero set of $\,x\cdot u\,$. The latter, a smooth
hypersurface in $\,Q\,$ which we call an
\textit{equator}, forms the  fixed-point set of
$\,\rho\,$, and separates $\,Q\,$ into two open topological
discs that we call \textit{hemispheres}.

Though a general ovaloid $\,Q\subset{\R^{n+1}}\,$ has no such
symmetry,  \textit{strict} ovaloids enjoy a very similar
involution relative to {any} hyperplane $\,P$ in
$\,{\R^{n+1}}$. Just as in the symmetric case, this
involution exchanges two hemispherical discs in $\,Q\,$ while
fixing the smooth ``equator'' that separates them. We can
define it very conveniently using Steiner
symmetrization:

\begin{definition}
\label{df:sym}  Suppose $\,K\subset{\R^{n+1}}\,$ is a compact
convex domain bounded by an ovaloid $\,Q$. Each line
perpendicular to a fixed hyperplane $P\subset{\R^{n+1}}\,$
intersects $\,K\,$ in a closed segment (possibly one point or
empty). Translating every such segment along the line
containing it until its midpoint lies on $\,P$, we map
$\,K\,$ to a new set which is clearly symmetric across
$\,P\,$. The boundary of this set again forms an ovaloid
\cite[Thm. 1.2.1]{gardner:GT}, which we label $\,Q_P\,$, and
call the \textit{Steiner $P$-symmetral} of $\,Q\,$. The
resulting map
\[
\sigma:Q\to Q_P\ ,
\]
is clearly continuous and injective, hence (by compactness of
$\,Q\,$) a homeomorphism. We call it the \textit{Steiner
$\,P$-symmetrization} of $\,Q\,$.
\end{definition}

We now use $\,\sigma\,$ to construct the promised involution
on $\,Q\,$.

\begin{definition}
\label{df:P-stuff}
The \textit{$P$-involution} on a strict $\,C^2\,$
ovaloid $\,Q\,$ is the map
\begin{equation}
\rho := \sigma^{-1}\circ\rp\circ\sigma\ ,
\label{eq:rho}
\end{equation}
where $\,\rp\,$ denotes reflection across $\,P\,$ (Eq.
\ref{eq:r1}). We call the fixed-point set
\[
\G:=
\left\{x\in Q\ :\ \rho(x) = x\ \right\}
\]
the \textit{$P$-equator} of $\,Q$.
\end{definition}
\medskip

\begin{lma}
\label{lma:gauss}
Given any hyperplane $\,P\,$ through the origin in
$\,\R^{n+1}\,$, we have
\begin{enumerate}
\item
The $P$-equator $\,\G\,$ on a
strict $\,C^2\,$ ovaloid $\,Q\subset\R^{n+1}\,$ comprises
those points of $\,Q\,$ having unit normals in $\,P\,$:
\[
\G = \nu^{-1}\left(P\cap\S^{n}\right)\ .
\]
\item
$\,\G\,$ is $\,C^1\,$ diffeomorphic to
$\,\S^{n-1}\,$, and
\item
$\,\G\,$ splits its complement in
$\,Q\,$ into two topological hemispheres, each forming the
graph of a $\,C^2\,$ function over a common open set in
$\,P\,$.
\end{enumerate}
\end{lma}

\begin{proof}
A line perpendicular to $\,P\,$ can intersect a strict ovaloid
$\,Q\,$ in only two ways: Either it (a) grazes $\,Q\,$
tangentially, in which case our involution $\,\rho\,$ fixes
the point of contact, or (b) it pierces $\,Q\,$ transversally
at exactly two points, and $\,\rho\,$ swaps this pair.  By
definition, the $P$-equator $\,\G\subset Q\,$ comprises the
points of case (a). So when $\,q\in\G\,$,   $\,T_qQ\,$
contains a line perpendicular to $\,P\,$, and we have
$\,\nu(q)\in P\,$. This proves the Lemma's first statment,
and since the unit normal map on a strict $\,C^2\,$ ovaloid
is a $\,C^1\,$ diffeomorphism, the second statement
follows easily.

Moreover, we see that each of the hemispherical regions
complementing $\,\G\,$ in $\,Q\,$ now comprises points
$\,q\,$ with $\,\nu(q)\not\in P\,$. By the Inverse Function
Theorem, this ensures that the orthogonal projection
$\,\pi:\R^{n+1}\to P\,$ restricts to a local $\,C^2\,$
diffeomorphism on each hemisphere. Since $\,\pi\,$ also
injects on each hemisphere into $\,P\,$, it now follows that
both hemispheres are $\,C^2\,$ graphs over the region bounded
by $\,\pi(\G)\subset\,P\,$. This verifies conclusion (3).
\end{proof}

We need one additional characterization of $\,\rho\,$
and $\,\G\,$.

\begin{lma}
\label{lma:rho}
We can express $\,\rho\,$ using a continuous ``signed
$\,P$-height'' function $\,\lambda:Q\to\R\,$ via
  \begin{equation}
  \rho(x) = x - 2\,\lambda(x)\, u \label{eq:lambda}\ .
  \end{equation}
Here $\,u\,$ denotes a unit vector normal
to $\,P\,$, and we have $\,\G = \lambda^{-1}(0)\,$.
\end{lma}

\begin{proof}
Using $\,\sigma\,$ to pull back Eq. (\ref{eq:r1})
(with $\rho\,$ replaced by $\,\rp\,$ there) from
$\,Q_P\,$, we get
\[
\rho(x) =  \sigma^{-1}\biggl(\sigma(x)-2(\sigma(x)\cdot
u)\, u\biggr)
\]
Since $\,\sigma\,$ acts by simple translation on the $\,
u$-parallel line through $\,x$, this immediately gives
(\ref{eq:lambda}), with
\[
\lambda(x) = \sigma(x)\cdot u
\]
Continuity of $\,\sigma\,$
makes $\,\lambda\,$ continuous too.
\end{proof}

\section{Main Results.}
\label{sec:main}

We now combine the facts we have developed about
lipschitz chains and ovaloids to produce a basic technical
result that we need to prove our main theorem.

\begin{prop}
\label{prop:OR}
Suppose we inject a compact, oriented riemannian manifold
$\,M^{n-1}\,$ into a strict $\,C^2\,$ ovaloid
$\,Q\subset\R^{n+1}\,$ with a lipschitz map $\,F\,$. If
$\,\pi_{\#}[F] =0\,$, then $\,\rho\,$ maps $\,F(M)\,$ to
itself, and reverses orientation; that is,
$\,F^{-1}\circ\rho\circ F\,$ constitutes an
orientation-reversing homeomorphism of $\,M\,$.
\end{prop}

\begin{proof}

We first note that $\,F(M)\,$ cannot lie wholly in
the $P$-equator $\,\G\,$. If it did, the injectivity of
$\,\pi\,$ on $\,\G\,$ would make $\,\pi\circ F\,$ injective
on $\,M\,$, and imply, by Lemma \ref{lma:spt}, that
$\,\spt(\pi_\#[F]) = \pi(F(M))\ne\emptyset\,$. This violates
our assumption that $\,\pi_\#[F]=0\,$.

We now argue that $\,\rho\,$ leaves $\,F(M)\,$
setwise fixed, as the Proposition claims. Indeed,
suppose not. Then since $\,F(M)\,$ doesn't lie wholly in
$\,\G\,$, we must have some point $\,q^+\,$ in one of the
hemispheres $\,\Omega^\pm\,$ complementing  $\,\G\,$ in
$\,Q\,$ for which
\begin{equation}
\label{eq:q+}
q^+\in F(M)\ ,\quad\text{but}\quad
q^- :=\rho(q^+)\not\in F(M)\ .
\end{equation}
Without loss of generality, assume $\,q^+\in\Omega^+\,$. Then
$\,q^-\,$ lies in $\,\Omega^-\,$, but, missing the compact set
$\,F( M)\,$, it must lie in an open set $\,U^-\subset
\Omega^-\,$ likewise disjoint from $\,F( M)\,$.

Define $\,U^+ := \rho(U^-)\subset \Omega^+\,$. Between the
third conclusion of Lemma~\ref{lma:gauss} and our definition
of $\,\rho\,$, we then see that both $\,U^+\,$ and $\,U^-\,$
are graphs of $\,C^2\,$ functions over a common open subset
of $\,U\subset P\,$. In other words,
$$
\pi(U^+) = \pi(U^-)=:U\subset P\ .
$$
Now observe that any 1-form $\,\phi\,$ supported in $\,U\,$
pulls back to $\,Q\,$, via $\,\pi\,$, as the sum of 1-forms
$\,\phi^+\,$ and $\,\phi^-\,$ supported in $\,U^+\,$ and
$\,U^-\,$ respectively:
\begin{equation*}
\pi^*\phi = \phi^++\phi^-\ ,\quad\text{with}\quad
\phi^+ = \left(\pi\big|_{U^+}\right)^*\phi
\quad\text{and}\quad
\phi^- = \left(\pi\big|_{U^-}\right)^*\phi\ .
\end{equation*}
Since $\,\pi|_{U^+}\,$ is a diffeomorphism, the relation
$\,\phi\leftrightarrow\phi^+\,$ induces an isomorphism
$\,\D^1(U)\approx\D^1(U^+)\,$. In view of Eq. (\ref{eq:q+})
above, Lemma~\ref{lma:spt} shows that the support of
$\,[F]\,$ contains $\,q^+\in U^+\,$. So there exists a 1-form
$\,\phi\,$ supported in $\,U\,$ such that $\,[F](\phi^+)\ne
0\,$. On the other hand, $\,[F](\phi^-)\,$ \emph{does}
vanish, because the support of $\,\phi^-\,$ lies in
$\,U^-\,$, which, by construction, entirely misses
$\,\spt([F])\,$. But then
$$
\left(\pi_{\#}[F]\right)(\phi)\ =\ [F](\pi^*\phi)\ =\
[F](\phi^++\phi^-)\ =\ [F](\phi^+)\ \ne\ 0\ .
$$
This contradicts the vanishing of
$\,\pi_{\#}[F]\,$, and thereby confirms our lemma's first
conclusion: $\,\rho\,$ preserves $\,F(M)\,$.

It also shows that the formula
\[
\Phi:=F^{-1}\circ\rho\circ F
\]
constructs a well-defined mapping of $\,M\,$. Since
both $\,F\,$ and $\,\rho\,$ are injective, and $\,M\,$ is
compact, $\,\Phi\,$ is bicontinuous. It remains to show that
it reverses orientation.

We started by ruling out the inclusion $\,F(M)\subset\G\,$.
But $\,F(M)\,$ cannot avoid $\,\G\,$ entirely, as this would
place it completely inside either $\,\Omega^+\,$, or
$\,\Omega^-\,$, where it could not be preserved by
$\,\rho\,$. Hence $\,F^{-1}(\G)\,$, which comprises the
fixed-point set of $\,\Phi\,$, is not empty. On the other hand
$\,\rho\,$ swaps $\,\Omega^+\,$ and $\,\Omega^-\,$,
so $\,\Phi\,$ moves every component of $\,M\setminus
F^{-1}(\G)\,$. The
brief note by Brown and Kister \cite{b&k:PAMS} now shows
that there are just two such components, which $\,\Phi\,$
must swap while reversing orientation.
\end{proof}
\bigskip

The main result discussed in our introduction,
Theorem~\ref{thm:main}, actually follows as a corollary of the
more complex statement below:

\begin{thm}
\label{thm:tech}
Consider a  a compact embedded lipschitz $(n-1)$-cycle
$\,F:M\to Q\,$ on a strict $\,C^2\,$ ovaloid
$\,Q\subset\R^{n+1}\,$. If $\,[F]\,$ projects to zero on
$\,n+1\,$ independent hyperplanes, then $\,M\,$ admits a
fixed-point-free homeomorphism of degree $\,(-1)^{n-1}\,$.
\end{thm}

\begin{proof}
We argue by contradiction: Suppose that for some compact,
embedded lipschitz $(n-1)$-cycle $\,F:M\to Q\,$, we had
$n+1\,$ independent hyperplanes $\,P_i\,$ whose corresponding
orthogonal projections $\,\pi_i:\R^{n+1}\to P_i\,$ all
annihilated $\,[F]$, i.e.,
$$
\pi_{i\#}[F]=0\qquad i=1,2,\dots,n+1
$$
Choose a unit normal vector $\, u_i\,$ for each $\,P_i$. Since
$\pi_i\,$ commutes with translation along $\,u_i\,$, we may
assume each hyperplane $\,P_i\,$ passes through the origin.
Lemma~\ref{lma:rho} now assigns a $\,P_i$-involution
$\,\rho_i\,$, along with a signed $P_i$-height function
$\,\lambda_i\,$, to $\,Q\,$ for each of these hyperplanes,
so that for each $\,i=1,\dots,n+1$, we have
\begin{equation}
\label{eqn:A}
\rho_i(x) = x - 2\,\lambda_i(x)\, u_i
\qquad\hbox{for all $\,x\in Q$}\ .
\end{equation}

Since $\,\pi_{i\#}[F] = 0$ for each $\,i\,$, we also know, by
Proposition~\ref{prop:OR}, that each of the homeomorphisms
$$
\Phi_i := F^{-1}\circ \rho_i\circ F
$$
reverses orientation on $\,M\,$. The composition
\begin{eqnarray*}
\Psi &:=& \Phi_{n+1}\circ \Phi_{n} \circ\cdots\circ
\Phi_2\circ \Phi_1\\
&=&F^{-1}\circ\rho_{n+1}\circ\rho_{n}\circ\cdots\circ\rho_2\circ\rho_1\circ F\ ,
\end{eqnarray*}
therefore preserves or reverses the orientation of $\,M\,$
depending respectively on whether $\,n\,$ is odd or even. In
other words, $\,\Psi:M\to M\,$ has degree $\,(-1)^{n+1} =
(-1)^{n-1}\,$, and we can complete the proof by showing that
$\,\Psi\,$ fixes no point in $\,M\,$.

Indeed, if we \emph{had} a fixed point $\,x\in M$, then
$\,y=F(x)\in Q\,$, would satisfy
\[
\rho_{n+1}\circ\rho_{n}\circ\cdots\circ\rho_2\circ\rho_1(y)
= y\ ,\quad\text{with}\quad y = F(x)\ .
\]
Expand this out using Eq. (\ref{eqn:A}), cancel the
lone $\,y$ on each side, and divide by $\,-2\,$, to get
\begin{eqnarray*}
\lefteqn{\lambda_1(y)\,u_1\ +\
\lambda_2\left(\rho_1(y)\right)\,u_2 \ +\
\lambda_3(\rho_2(\rho_1(y)))\ + \ \cdots}\\ & &\quad\cdots\
+\
\lambda_{n+1}(\rho_{n}(\rho_{n-1}(\cdots(\rho_1(y))\cdots)))\,u_{n+1}
= 0\ .
\end{eqnarray*}
Since the $\,u_i$'s are linearly independent, their
coefficients above must all vanish.  Working recursively from
left to right using the last conclusion of
Lemma~\ref{lma:rho} and the fact that $\,\rho_i\,$ fixes
$\,\G_i\,$, we now reason as follows:
$$
\lambda_1(y) = 0\quad\Longrightarrow\quad
y\in\Gamma_1\quad\hbox{and}\quad \rho_1(y)=y\ .
$$
Knowing now that $\,\rho_1(y)=y$, we subsequently get
$$
\lambda_{2}(y) = \lambda_2(\rho_1(y)) =
0\quad\Longrightarrow\quad y\in\Gamma_2\quad\hbox{and}\quad
\rho_2(y)=y\ .
$$
Continuing in this way, we find that for \emph{all} $\,i =
1,\,2,\,\dots,\,n+1\,$, we have  $\,\lambda_i(y)=0$ and
$\,y\in\Gamma_i\,$. It then follows from the first
conclusion of Lemma~\ref{lma:gauss} consequently that
$$
y\in \Gamma_1\cap\Gamma_2\cap\cdots\cap\Gamma_{n+1} =
\nu^{-1}\left(S^{n+1}\cap P_1\cap P_2\cap\cdots\cap
P_{n+1}\right)\ .
$$
But the intersection on the right is empty, because $\,n+1\,$
independent hyperplanes in $\,\R^{n+1}\,$ intersect only at
the origin. We have thus contradicted the existence of a
fixed-point for $\,\Psi\,$, and proven the theorem.
\end{proof}

Our main result now follows easily:

\begin{thm}
\label{thm:main}
A lipschitz embedding of a homology $(n-1)$-sphere on a strict
$\,C^2\,$ ovaloid in $\,\R^{n+1}\,$ cannot project to zero on
$\,n+1\,$ independent hyperplanes.
\end{thm}

\begin{proof}
By Theorem \ref{thm:tech}, a violation of this theorem would
produce a fixed-point-free homeomorphism of an
($n-1$)-dimensional homology sphere with degree
$\,(-1)^{n-1}\,$. By a well-known consequence of the
Lefschetz fixed-point theorem \cite[Cor. 6.21]{vick:HT}, no
such map exists.
\end{proof}

\section*{Acknowledgements}
The author thanks Mohammad Ghomi, William P. Ziemer, and Chuck
Livingston for conversations pertinent to this paper; Allan
Edmonds,
for pointing out reference \cite{b&k:PAMS} (which simplified a
key point in the proof of Theorem \ref{thm:tech}) and Stanford
University for supporting part of this project.

\providecommand{\bysame}{\leavevmode\hbox to3em{\hrulefill}\thinspace}

\vfill

\vskip .5in

\def\x{{\vrule height 5pt width 9pt depth 3pt}}
\def\o{{\vrule height 0pt width 9pt depth 0pt}}
\baselineskip 0pt
\centerline{\o\x\x\x\x\o\x\x}
\centerline{\o\x\o\x\x\o\x\o}
\centerline{\x\x\o\x\x\x\x\o}
\medskip

\end{document}